\documentclass{article}
\usepackage{amsmath,amssymb,amsthm}
\newcommand{\var}{\text{\rm var }}
\newcommand{\E}{\mathbb E}

\newcommand{\R}{\mathbb R}

\newcommand{\x}{\mathbf x}

\newcommand{\dd}{\text{\rm d}}

\newtheorem{theorem}{Theorem}
\newtheorem{lemma}{Lemma}

\newtheorem{example}{Example}
\newtheorem{remark}{Remark}
\newtheorem{proposition}{Proposition}
\begin{document}
\title{Interaction processes for unions of facets,\\ a limit behavior}
\author{Jakub Ve\v{c}e\v{r}a, Viktor Bene\v{s}}

\maketitle
\noindent Charles University in Prague, Faculty of Mathematics and Physics, Department of Probability and Mathematical Statistics, Sokolovsk\'{a} 83, 18675 Praha 8-Karl\'{\i}n, Czech Republic, vecera@karlin.mff.cuni.cz, benesv@karlin.mff.cuni.cz

\bigskip

\noindent{\bf Abstract}

\medskip

\noindent 
In the series of models with interacting particles in stochastic geometry, a further contribution presents the facet process which is defined in arbitrary Euclidean dimension. In 2D, 3D specially it is a process of interacting segments, flat surfaces, respectively. Its investigation is based on the theory of functionals of finite spatial point processes given by a density with respect to a Poisson process. The methodology based on $L_2$ expansion of the covariance of functionals of Poisson process is developed for $U$-statistics of facet intersections which are building blocks of the model. The importance of the concept of correlation functions of arbitrary order is emphasized. Some basic properties of facet processes, such as local stability and repulsivness are shown and a standard simulation algorithm mentioned. Further the situation when the intensity of the process tends to infinity is studied. In the case of Poisson processes a central limit theorem follows from recent results of Wiener-Ito theory. In the case of non-Poisson processes we restrict to models with finitely many orientations. Detailed analysis of correlation functions exhibits various asymptotics for different combination of $U$-statistics and submodels of the facet process.

\medskip


\medskip


\section{Introduction}
Modeling of particle systems with interactions by means of finite point processes presents an interesting field of stochastic geometry the research of which is not yet completed. The work of \cite{RefK} concerning area-interaction models of planar discs was later developed by \cite{RefM} to models of interacting discs with more general densities (with respect to the Poisson process) from exponential families. The authors also developed sophisticated simulation procedures of the models. The motivation for the modeling was to serve for statistical purposes in the evaluation of real data.

In the present paper we develop similar interaction processes in a bounded window in arbitrary Euclidean dimension. The difference in comparison to \cite{RefM} is that the particles (called facets) are lower dimensional, they form compact subsets of hyperplanes. The interactions arise from intersections of facets and the their global amounts form $U$-statistics \cite{RefR}.
In 2D and 3D these models may also serve to real data evaluation of segment, surface processes, respectively. Rather than statistics we develop the facet process theory in arbitrary dimension and deal with limit behavior when the intensity of the reference Poisson process tends to infinity.

Recently in \cite{RefB} functionals of spatial point processes given by a density with respect to the Poisson process were investigated using the $L_2$ expansion from \cite{RefL} which is applied to the product of a functional and the density. Using a special class of functionals called $U$-statistics closed formulas for mixed moments of functionals are obtained. Similar formulas, but under a stronger assumption of a product form of the driving function of the functional, were derived in \cite{RefD} using the Georgii-Nguyen-Zessin formula. In processes with densities the key characteristics are the correlation functions \cite{RefG} of arbitrary order which are dual to kernel functions of the density as a function of the Poisson process. 
 
We call facets some bounded subsets of hyperplanes with a given shape, size and orientation. Natural geometrical characteristics of the union of the facets, based on Hausdorff measure of the intersections of pairs, triplets, etc., of facets form $U$-statistics. Building a parametric density from exponential family, the limitations for the space of parameters have to be given. In the paper basic properties of facet processes are investigated, based on the fact that the densities come from exponential family as studied in \cite{RefM} for the disc process with an analogous density type. It is shown that the given class of processes is characterized by a repulsive behavior. Conditions for local stability are given and a standard simulation algorithm based on birth-death algorithm in Markov chain Monte Carlo is mentioned, which enables to observe realizations of the facet process.

As an application of the moment formulas we are interested in the limit behaviour when the intensity of the reference Poisson process tends to infinity. Central limit theorems for $U$-statistics of Poisson processes were derived based on Malliavin calculus and the Stein method in \cite{RefR}. The results were extended to the multivariate case in \cite{RefL} and they apply to the vector of model characteristics of the Poisson facet process in our setting. In the paper \cite{RefX} the authors give conditions for functionals of Gibbs point processes under which they are asymptotically Gaussian with increasing window. Another related work is \cite{RefP}, where the central limit theorem for the number of intersections in expanding window is derived for the stationary planar segment process satisfying certain conditions on absolute regularity coefficients.

In facet processes having densities we restrict ourselves to the model with finitely many orientations corresponding to canonical vectors. Submodels of the facet process are investigated where a detailed analysis of the correlation function yields different asymptotics. When the order of the submodel is not greater than the order of the observed $U$-statistic then asymptoticaly the mean value of the $U$-statistic vanishes. This leads to a degeneracy in the sense that some orientations are missing. On the other hand  when the order of the submodel is greater than the order of the observed $U$-statistic then the limit of correlation function is finite and nonzero and under selected standardization $U$-statistic tends almost surely to its non-zero expectation. By changing the standardization, however, we achieve a finite non-zero asymptotic variance. Even if these results are obtained in a special situation with facets of a fixed shape and size related to the window size, it is important that they allow us to understand the ongoing problems for a possible further investigation of the complex model.

\section{The background of point processes having a density}\label{sec:1} 
Consider a bounded Borel set $B\subset {\mathbb R}^d$ with Lebesgue measure $|B|>0$ and a measurable space $({\mathbf N},{\mathcal N})$ of integer-valued finite measures on $B.$ ${\mathcal N}$ is the smallest $\sigma $-algebra which makes the mappings $\x\mapsto \x(A)$ measurable for all Borel sets $A\subset B$ and all $\x\in {\mathbf N}.$ A random element having a.s. values in $({\mathbf N},{\mathcal N})$ is called a finite point process. Integer-valued finite measures can be represented by systems of points corresponding to their support. Let a Poisson point process $\eta $ on $B$ have finite intensity measure $\lambda $ with no atoms, $\lambda (B)>0,$ and distribution $P_\eta $ on $\mathcal N.$ For a measurable map $F:{\mathbf N}\rightarrow {\mathbb R}$ it holds \cite{RefA}
 \begin{equation} \label{odhad} {\mathbb E}[F(\eta )]=e^{-\lambda (B)}\sum_{n=0}^\infty\frac{1}{n!}\int_B\dots\int_B F(u_1,\dots ,u_n)\lambda^n (\dd(u_1,\dots ,u_n)), \end{equation} where we write $\lambda^n (\dd (x_1,\dots ,x_n))$ instead of $\lambda (\dd x_1)\dots\lambda (\dd x_n).$ 
Further we consider a finite point process $\mu $ on $B$ given by a density $p$ w.r.t. $\eta ,$ i.e. with distribution $P_\mu $ \begin{equation}\label{dns}dP_\mu (\x)=p(\x)dP_\eta (\x),\; \x\in \mathbf N,\end{equation} where $p:{\mathbf N}\rightarrow {\mathbb R}_+$ is measurable satisfying $$\int_{\mathbf N}p(\x)\dd P_\eta (\x)=1.$$ The consequence of (\ref{dns}) is a formula \begin{equation}\label{odha}{\mathbb E}F(\mu )={\mathbb E}[F(\eta )p(\eta )].\end{equation}
Let $\mu $ be a finite point process with density $p$ satisfying \begin{equation}\label{hered}p(\x)>0\Rightarrow p(\tilde{\x})>0\end{equation} for all $\tilde{\x}\subset \x.$ For the (Papangelou) conditional intensity of $\mu ,$ see \cite{RefA}, it holds $$\lambda^*(u;\x)=\frac{p(\x\cup \{u\})}{p(\x)},\; \x\in {\mathbf N},\: u\in B,\: u\notin \x.$$ For $p(\x)=0$ we put $\lambda^*(u;\x)=0.$ For $n>1$ we use analogously $$\lambda^*_n(u_1,\dots ,u_n;\x)=\frac{p(\x\cup \{u_1,\dots ,u_n\})}{p(\x)},$$ $u_1,\dots ,u_n\in B\setminus\{\x\}
$ distinct, the conditional intensity of $n$-th order of $\mu ,\;\lambda^*_0\equiv 1.$ We observe that $\lambda^*_n$ is symmetric in the variables $u_1,\dots ,u_n.$ The expectations of conditional intensities \begin{equation}\label{coref}\rho_n(u_1,\dots ,u_n;\mu )={\mathbb E}\lambda^*(u_1,\dots ,u_n;\mu )={\mathbb E}[p(\eta\cup\{u_1,\dots ,u_n\})],\end{equation} are called $n-$th correlation functions of the point process $\mu ,$ cf. \cite{RefG}.

For a functional $F,\;y\in B,$ one defines the difference operator $D_yF$ for a point process $\mu $ as a random variable $$D_yF(\mu )=F(\mu + \delta_y)-F(\mu ),$$ where $\delta_y $ is a Dirac measure at the point $y.$ Inductively for $n\geq 2$ and $(y_1,\dots ,y_n)\in B^n$ we define a function $$D_{y_1,\dots ,y_n}^nF=D_{y_1}^1D_{y_2,\dots ,y_n}^{n-1}F,$$
where $D^1_y=D_y,\; D^0F=F.$ Operator $D^n_{y_1,\dots ,y_n}$ is symmetric in $y_1,\dots ,y_n$ and symmetric functions $T_n^\mu F$ on $B^n$ are defined as $$T^\mu_nF(y_1,\dots ,y_n)={\mathbb E}D_{y_1,\dots ,y_n}^nF(\mu ),$$ $n\in {\mathbb N},\; T_0^\mu F={\mathbb E}F(\mu ),$ whenever the expectations exist. We write $T_nF$ for $T^\eta_nF.$ For the functionals of a Poisson process Theorem 1.1 in \cite{RefJ} says that given $F,\tilde{F}\in L^2(P_\eta )$ it holds
\begin{equation}\label{skaso}{\mathbb E}[F(\eta ){\tilde F}(\eta )]={\mathbb E}F(\eta ){\mathbb E}{\tilde F}(\eta )+\sum_{n=1}^\infty\frac{1}{n!}\langle T_nF,T_n{\tilde F}\rangle_n,\end{equation} where $\langle .,.\rangle_n$ is the scalar product in $L_2(\lambda^n).$ We will use symbol $[n]=\{1,\dots ,n\},$ $n\in {\mathbb N}.$

For  $p \in L_2(P_\eta ),\;n\in {\mathbb N},$ it holds \begin{equation}\label{tnp}T_np(y_1,\dots , y_n)=\sum_{J\subset [n]}(-1)^{n-|J|}\rho_{|J|}(\{y_j,\:j\in J\};\mu )\end{equation} for $\lambda^n$-almost all $(y_1,\dots , y_n),$ where $|J|$ is the cardinality of $J.$

\section{U-statistics} 
A $U$-statistic of order $k\in {\mathbb N}$ of a finite point process $\mu $ is a functional defined by
\begin{equation}\label{ust}F(\mu )=\sum_{(x_1,\dots ,x_k)\in\mu^k_{\neq }} f(x_1,\dots ,x_k),\end{equation} where $f:B^k\rightarrow {\mathbb R}$ is a function symmetric w.r.t. to the permutations of its variables, $f\in L_1(\lambda^k ).$ Here $\mu^k_{\neq }$ is the set of $k$-tuples of different points of $\mu .$ We say that $F$ is driven by $f.$ In this Section basic results on $U$-statistics for point processes having densities, obtained in \cite{RefB} are reviewed for later use. By the Slivnyak-Mecke theorem \cite{RefS} we have 
$${\mathbb E}F(\eta )=\int_{B^k}f(x_1,\dots ,x_k)\lambda^k (\dd (x_1,\dots ,x_k)).$$ 
For a $U$-statistic $F\in L_2(P_\eta )$ of order $k$ and density $p\in L_2(P_\eta )$ it holds \begin{equation}\label{effa}{\mathbb E}F(\mu )=\int_{B^k}f(x_1,\dots ,x_k)\rho_k(x_1,\dots ,x_k;\mu)\lambda^k (\dd (x_1,\dots ,x_k)).\end{equation} In integrals like (\ref{effa}) it does not mind that $\rho_k$ is defined only for distinct arguments $x_j.$  We are interested in higher-order and mixed moments in the following.

 We can use a short expression of formulas for moments using diagrams and partitions, see \cite{RefT}, \cite{RefL}. Let $\tilde{\prod}_k$ be the set of all  partitions $\{J_i\}$ of $[k],$ where $J_i$ are disjoint blocks and $\cup J_i=[k].$ For
$k=k_1+\dots +k_m$ and blocks $$J_i=\{ j: k_1+\dots +k_{i-1}< j\leq k_1+\dots +k_i\},\; i=1,\dots ,m,$$ consider the partition $\pi =\{J_i,\;1\leq i\leq m\}$ and let
$\prod_{k_1,\dots ,k_m}\subset \tilde{\prod}_k$ be the set of all partitions $\sigma\in\tilde{\prod}_k$ such that $|J\cap J'|\leq 1$ for all $J\in\pi $ and all $J'\in\sigma .$ Here $|J|$ is the cardinality of a block $J\in\sigma .$ 
For a partition $\sigma\in\prod_{k_1\dots k_m}$ we define the function $(\otimes_{j=1}^mf_j)_\sigma:B^{|\sigma |}\rightarrow {\mathbb R}$ by replacing all variables of the tensor product $\otimes_{j=1}^mf_j$ that belong to the same block of $\sigma $ by a new common variable, $|\sigma |$ is the number of blocks in $\sigma .$
\begin{theorem}\label{theo2}Let $m\in\mathbb{N},\;\prod_{i=1}^mF_i\in L_2(P_{\eta}),$ $p\in L_2(P_\eta),$ where $F_i$ are $U$-statistics of orders $k_i$ driven by nonnegative functions $f_i,$ respectively, $i=1,\dots,m$. Then \begin{equation}\label{divi}{\mathbb E}\left[\prod_{i=1}^mF_i(\mu )\right]=\sum_{\sigma\in\prod_{k_1\dots k_m}}\int_{B^{|\sigma |}}(\otimes_{i=1}^mf_i)_{\sigma }(x_1,\dots ,x_{|\sigma |})\times \end{equation}$$\times \rho_{|\sigma |}(x_1,\dots ,x_{|\sigma |};\mu )\lambda(\dd (x_1,\dots ,x_{|\sigma |})).$$\end{theorem}

Similar formulas to (\ref{divi}), but under a stronger assumption of a product form of the driving function $f$ of $F,$ were derived in \cite{RefD} using the Georgii-Nguyen-Zessin formula \cite{RefA}. As an application consider processes $\mu_a$ with densities $p_a$ w.r.t. $\eta_a,\,a\geq 1.$  Let formula (\ref{divi}) be applied with $ \lambda_a.$ The term at the highest power of $a$ which comes from $a\lambda $ is called the leading term. The rate of the leading term is $\sum_{i=1}^{m}k_{i}$. There remains the dependence on $a$ hidden in $\rho ,$ which will be investigated later. \begin{theorem}\label{tht}
For a $U$-statistic $$F(\mu_a)=\sum_{(x_1,\dots ,x_{k})\in\mu^{k}_{a\neq }} f(x_1,\dots ,x_{k})$$ of order $k$ on a bounded set $B$ we have the following expression for the leading term of centered moments, the rate being $mk:$
 ${\mathbb E}[(F(\mu_a)-{\mathbb E}F(\mu_a))^m]=$ $$\sum_{l=0}^m\binom{m}{l} (-1)^{m-l}\int_{B^{lk}} f^{\otimes l}(x_1,\dots ,x_{lk})\rho_{lk}(x_1,\dots ,x_{lk};\mu_a)\lambda^{lk}(\dd (x_1,\dots ,x_{lk})) $$$$\times\left(\int_{B^k}f(x_1,\dots ,x_{k})\rho_{k}(x_1,\dots ,x_{k};\mu_a)\lambda^{k}(\dd (x_1,\dots ,x_{k}))\right)^{m-l}.$$ \end{theorem}
\noindent{\bf Proof:} We have $${\mathbb E}[(F(\mu_a)-{\mathbb E} F(\mu_a))^m]=\sum_{l=0}^m\binom{m}{l} (-1)^{m-l}{\mathbb E}F(\mu_a)^l({\mathbb E} F(\mu_a))^{m-l}.$$
From (\ref{effa}) it is $({\mathbb E} F(\mu_a))^{m-l}=$$$=a^{k(m-l)}\left(\int_{B^k}f(x_1,\dots ,x_{k})\rho_{k}(x_1,\dots ,x_{k};\mu_a)\lambda^{k}(\dd (x_1,\dots ,x_{k}))\right)^{m-l}$$ and from (\ref{divi}) applied to ${\mathbb E}[F(\mu_a)^l]$ we take the term with highest
power of $a.$ It corresponds to $$a^{lk}\int_{B^{lk}} f^{\otimes l}(x_1,\dots ,x_{lk})\rho_{lk}(x_1,\dots ,x_{lk};\mu_a)\lambda^{lk}(\dd (x_1,\dots ,x_{lk})).$$ These terms come from $\sigma\in\prod_{=1,k,\dots ,k}\subset\prod_{k,\dots ,k} ,$ where $\prod_{=1,k,\dots ,k}$ is a singleton with the only partition containing $mk$ blocks all of cardinality one.\hfill $\Box $

Let $B\subset {\mathbb R}^d,\,d\in {\mathbb N}$ be as in Section 2, $X$ be a point process of compact sets \cite{RefS} called grains and we denote $\x$ realization of $X$ on $B$, i.e. the union of all compact sets. Consider that there is a probability density \cite{RefM} \begin{equation}\label{den}p(\x)=c_\nu\exp(\nu\cdot G(\x)), \end{equation}of $X$
w.r.t. a given reference Poisson point process $\eta $ of compact sets. Here $\nu =(\nu_1,\dots ,\nu_d)$ is a vector of real parameters, $c_\nu $ is a normalizing constant, $$G(\x)=(G_1({\x}),\dots,G_d({\x}))\in {\mathbb R}^d$$ is a vector of geometrical characteristics of $\x.$ In the exponent of (\ref{den}) there is the inner product in ${\mathbb R}^d.$
The largest set of ${\nu}$ such that exponential family density (\ref{den}) is well defined is
$$\{\nu\in\mathbb{R}^d : \mathbb{E}[\exp(\nu\cdot G(\eta ))] <\infty\}.$$
For  $m\in\mathbb{N}$
denote $$D^m_{y_1,\dots,y_m}G(\x)=(D^m_{y_1,\dots,y_m}G_1({\x}),\dots,D^m_{y_1,\dots,y_m}G_d({\x}))^T$$
the vector of $m-$th differences.
\begin{proposition}\label{Th:papan}Consider the probability density (\ref{den}).
 Then for the corresponding conditional
intensity $\lambda^*_m$ of order $m\in {\mathbb N}$ and a realization $\x$ of $X$ it holds
\begin{equation}
\lambda^*_m(y_1,\dots,y_m;\x)=\exp (\nu\cdot Q_mG(\x)),
\end{equation}
where $y_1,\dots ,y_m\in B\setminus\x$ are distinct,
\begin{eqnarray*}
Q_mG(\x)&=&D^m_{y_1,\dots,y_m}G(\x)\\
&&+\sum_{i_1,\dots,i_{m-1}\in\{1,\dots,m\}}D^{m-1}_{y_{i_1},\dots,y_{i_{m-1}}}G(\x)+\dots+
\sum_{1\leq i\leq m}D_{y_i}G(\x).
\end{eqnarray*}
\end{proposition}

The intensity of the reference Poisson process depends on a specific model, see \cite{RefM} for interacting discs.
\section{Facet processes with density}
 Here we consider processes of interacting facets in ${\mathbb R}^d,\,d\in {\mathbb N},$ with their natural $U$-statistics. Let \begin{equation}\label{bck}Y=B\times (0,b]\times {\mathbb S}^{d-1},\end{equation} where $b>0$ is a size parameter, ${\mathbb S}^{d-1}$ is the hemisphere of axial orientations in ${\mathbb R}^d.$ For a point $y\in Y,\;y=(z,r,\phi )$ represents a facet which is a subset of a hyperplane $A=A(z,\phi)$ through point $z$ having normal orientation $\phi .$ Then \begin{equation}\label{dist}y=\{s\in A;\, dist(s,z)\leq r\}\end{equation} for a distance $dist$ in $A.$ The Poisson process $\eta $ on $Y$ has intensity measure $\lambda ,$ \begin{equation}\label{inti}\lambda (\dd (z, r, \phi ))=\chi (z)\dd zQ(\dd r)v(\phi )(\dd\phi ),\end{equation} where $Q$ is the size distribution, a probability measure on $(0,b],$ $v$ is a probability density w.r.t. Lebesgue measure on ${\mathbb S}^{d-1},$ $\chi $ is a bounded intensity function of facet centres on $B.$ 

Let $k\in\{2,\dots ,d\}, $ we consider the intersection of any $k$ facets such that the corresponding hyperplanes are in general position, cf. \cite{RefS}, p.133. The $\lambda^k$-measure of all such $k$-tuples is equal to $\lambda (Y)^k$ for each $k$ since we have a density $v$ in (\ref{inti}). For the Hausdorff measure ${\mathbb H}^j$ in ${\mathbb R}^d$ of order $j$ and $\x\in{\mathbf N}$ (on $Y$) we put \newpage\begin{equation}\label{hrm}G_1(\x )=\sum_{y\in \x }{\mathbb H}^{d-1}(y),\end{equation}
$$G_2(\x )=\frac{1}{2}\sum_{(y_1,y_2)\in \x^2_{\neq } }{\mathbb H}^{d-2}(y_1\cap y_2),$$ \centerline\vdots
$$G_{d-1}(\x )=\frac{1}{(d-1)!}\sum_{(y_1,\dots ,y_{d-1})\in \x^{d-1}_{\neq } }{\mathbb H}^1(\cap_{i=1}^{d-1}y_i),$$
$$G_{d}(\x )=\frac{1}{d!}\sum_{(y_1,\dots ,y_{d})\in \x^{d}_{\neq } }{\mathbb H}^0(\cap_{i=1}^{d}y_i).$$
We deal with \begin{equation}\label{gseg}G(\x)=(G_1(\x),\dots , G_d(\x)),\end{equation} here $G_j$  is $U$-statistic of the $j-$th order, $j=1,\dots ,d.$ The facet process $\mu $ is defined by the density (\ref{den}) $$p(\x)=c_\nu\exp(\nu\cdot G(\x))$$ with respect to $\eta ,$ where $\nu = (\nu_1,\dots ,\nu_d)\in\R^d.$  For $\x=\{u_1,\dots ,u_n\}$ we have $$G_j(\x)\leq const.n^jb^{d-j},$$ which together with the assumption $\nu_j\leq 0,\,j=2,\dots ,d$  proves using (\ref{odhad}) that $$p\in L_1(P_\eta )\cap L_2(P_\eta ).$$ In fact the assumption that the orientation distribution $V(d\phi )$ has a density $v$ can be weakened to the assumption that $V$ has at least $d$ atoms, which we will use later. \begin{example}Specially for $d=3$ the facet process may serve as a model for platelike particles in materials microstructure of metals. Here $G_1$ yields the total
 area of all plates, $G_2$ is the total length of intersection segments of pairs of particles, $G_3$ is the total number of intersections of triplets of particles. The size of negative parameter $\nu_2$ or $\nu_3$ gives the measure of neglection of intersections (repulsion) of the corresponding type. If $\nu $ is a zero vector, then the facet process is Poisson (with no repulsion).
\end{example}
We will consider special types of facet processes $\mu^{(k)},\; k=1,\dots ,d,$ such that in (\ref{den}) we have $\nu_k\leq 0$ while $\nu_j=0,\; j\neq k.$ We say that $\mu^{(k)}$ is a submodel of order $k$, especially  $\mu^{(1)}$ is a Poisson process and 
 $\mu^{(k)}$ is a Poisson process if and only if $\nu_k=0,\;k=2,\dots ,d.$ 
 
Local stability of a facet process is defined by the existence of $\alpha >0$ such that 
\begin{gather}\label{locst}\lambda^*(u;\x)\leq\alpha ,\; \x\in{\bf N},\;u\in Y\setminus\{\x\} 
\end{gather}
 and under given conditions the processes $\mu^{(k)},\; k=1,\dots ,d$ are locally stable. Simulation of a facet process $\mu $ is available using the birth-death Metropolis-Hastings algorithm \cite{RefGM}.

In \cite{RefM} the process with conditional intensity $\lambda^*$ is called attractive if $$\lambda^*(y;\x_1)\geq \lambda^*(y;\x_2),$$ for any
$\x_1,\x_2\in {\bf N},\;\x_2\subset \x_1.$ It is repulsive if $$\lambda^*(y;\x_1)\leq \lambda^*(y;\x_2),$$ for any
$\x_1,\x_2\in {\bf N},\;\x_2\subset \x_1.$ In the case of sharp inequality we say that the process is strictly attractive, strictly repulsive, respectively. 
\begin{proposition} The facet process $\mu^{(1)}$ is neither strictly attractive nor strictly repulsive for any $\nu_1\in\R.$
The facet processes $\mu^{(k)}$ are repulsive for $k=2,\dots ,d.$\end{proposition}

\section{Asymptotics with increasing intensity}
\subsection{The Poisson case}
 Generally on $B\subset {\mathbb R}^d,\,\lambda $ as in Section 2, for $l\geq 1$ and $i=1,\dots ,l$ let $k_i\in {\mathbb N},$ $f^{(i)}\in L_1(\lambda^{k_i})$ be symmetric functions.  Consider Poisson processes $\eta_a$ with intensity measures $\lambda_a=a\lambda ,\;a\geq 1.$ 
Following \cite{RefL} $U$-statistics \begin{equation}\label{poispos}F^{(i)}_a(\eta_a)=\sum_{(x_1,\dots ,x_{k_i})\in\eta^{k_i}_{a\neq }} f^{(i)}(x_1,\dots ,x_{k_i})\end{equation} are transformed to \begin{equation}\label{tri}
\hat{F}_a^{(i)}=a^{-(k_i-\frac{1}{2})}(F_a^{(i)}-{\mathbb E}F_a^{(i)}).\end{equation}
The asymptotic covariances are
\begin{equation}\label{asco}C_{ij}=\lim_{a\rightarrow\infty }cov(\hat{F}_a^{(i)},\hat{F}_a^{(j)})=\int_B T_1F^{(i)}(x)T_1F^{(j)}(x) \lambda (\dd x),\;i,j\in\{1,\dots ,l\}.\end{equation} 
The convergence under the distance between 
 $l$-dimensional random vectors $U,Z$ 
$$d_3 (U,Z)=\sup_{g\in {\cal H}}|{\mathbb E}g(U)-{\mathbb E}g(Z)|,$$ where $\cal H$ is the system of functions $h\in C^3({\mathbb R}^l)$ with 
\begin{gather}
\nonumber \max_{1\leq i_1\leq i_2\leq l}\sup_{x\in {\mathbb R}^l}\big|\frac{\partial^2h(x)}{\partial x_{i_1}\partial x_{i_2}}\big|\leq 1,\quad
\max_{1\leq i_1\leq i_2\leq i_3\leq l}\sup_{x\in {\mathbb R}^l}|\frac{\partial^3h(x)}{\partial x_{i_1}\partial x_{i_2}\partial x_{i_3}}|\leq 1
\end{gather}
 implies convergence in distribution. Based on the multi-dimensional Malliavin-Stein inequality for the distance $d_3$ of a random vector from a centered Gaussian random vector $Z$ with covariance matrix $C=(C_{ij})_{i,j=1,\dots ,l},$ 
\cite{RefL} show that under the assumption \begin{equation}\label{assd}\int_B |T_1F^{(i)}|^3\dd\lambda <\infty,\;i=1,\dots ,l,\end{equation} there exists a constant $c$ such that \begin{equation}\label{mltv}d_3((\hat{F}_a^{(1)},\dots ,\hat{F}_a^{(l)}),Z)\leq ca^{-\frac{1}{2}},\; a\geq 1.\end{equation}
\begin{example}For the Poisson facet processes $\eta_a,\,a\geq 1$ on $Y$ (\ref{bck}) with intensity measure $a\lambda $ (\ref{inti}) and the $U$-statistics $G_j(\eta), j=1,\dots ,d,$ in (\ref{hrm}) we obtain that $$T_1G_j(x)=\frac{1}{(j-1)!}\int_Y\dots\int_Y{\mathbb H}^{d-j}(\cap_{i=1}^{j-1}y_i\cap x)\lambda^{j-1} (\dd (y_1,\dots ,y_{j-1})).$$ The finiteness of the intensity measure $\lambda $ in (\ref{inti}) and the boundedness of the facets guarantee that all integrals (\ref{assd}) and (\ref{asco}) are finite. Thus for the random vector $(G_j(\eta_a),j=1,\dots ,d)$ both the central limit theorem when $a\rightarrow\infty $ and the Berry-Esseen type inequality (\ref{mltv}) hold. 
\end{example}
\subsection{The non-Poisson case}
Let facet processes $\mu_a,\,a\geq 1,$ have densities \begin{equation}
p_a(\x)=c_{\nu ,a}\exp (\nu\cdot G(\x)) 
\end{equation} w.r.t. Poisson processes $\eta_a$ with intensities $\lambda_a=a\lambda ,$ respectively. Here $G(\x)$ is given in (\ref{hrm}) and $\nu_j\leq 0,\; j=2,\dots ,d.$
We investigate $U$-statistics and submodels of the same order $d-k,$ from formula (\ref{effa}) we have
\begin{equation}\label{znovu}{\mathbb E} G_{d-k}(\mu_a^{(d-k)})=\end{equation}$$=\frac{a^{d-k}}{(d-k)!}\int_{Y^{d-k}}{\mathbb H}^k(\cap_{i=1}^{d-k}y_i)\rho_{d-k}(y_1,\dots ,y_{d-k};\mu_a^{(d-k)})\lambda^{d-k}(\dd (y_1,\dots ,y_{d-k})),$$$ k=0,\dots ,d-2.$ 
\begin{lemma}\label{ll1}For the processes $\mu_a^{(d-k)}$ the correlation function $\rho_{d-k},\;k=0,\dots ,$ $d-2,$ has form \begin{equation}\label{rods}
\rho_{d-k}(y_1,\dots ,y_{d-k};\mu_a^{(d-k)})=\frac{A(a)}{B(a)},\end{equation} where $A(a)=$
$$\sum_{n=0}^\infty \frac{a^n}{n!}\int_Y\dots \int_Y\exp\left(\frac{\nu_{d-k}}{(d-k)!}\sum_{(x_1,\dots ,x_{d-k})\subset\{u_1,\dots ,u_n,y_1,\dots ,y_{d-k}\}}{\mathbb H}^k(\cap_{i=1}^{d-k}x_i)\right)$$$$\times\lambda^n (\dd(u_1,\dots ,u_n)),$$$$B(a)=\sum_{n=0}^\infty \frac{a^n}{n!}\int_Y\dots \int_Y\exp\left(\frac{\nu_{d-k}}{(d-k)!}\sum_{(x_1,\dots ,x_{d-k})\subset\{u_1,\dots ,u_n\}}{\mathbb H}^k(\cap_{i=1}^{d-k}x_i)\right)$$$$\times\lambda^n (\dd(u_1,\dots ,u_n)).$$\end{lemma}
\noindent{\bf Proof:} Using formulas (\ref{coref}),(\ref{odhad}) we obtain
$$\rho_{d-k}(y_1,\dots ,y_{d-k};\mu_a^{(d-k)})=c_{\nu ,a}e^{-a\lambda (Y)}\sum_{n=0}^\infty \frac{a^n}{n!}\int_Y\dots \int_Y$$$$\exp\left(\frac{\nu_{d-k}}{(d-k)!}\sum_{(x_1,\dots ,x_{d-k})\subset\{u_1,\dots ,u_n,y_1,\dots ,y_{d-k}\}}{\mathbb H}^k(\cap_{i=1}^{d-k}x_i)\right)\lambda^n (\dd(u_1,\dots ,u_n)).$$

The normalizing constant can be expressed from
$$1={\mathbb E} p(\eta_a)=c_{\nu ,a}e^{-a\lambda (Y)}\sum_{n=0}^\infty \frac{a^n}{n!}\int_Y\dots \int_Y$$$$\exp\left(\frac{\nu_{d-k}}{(d-k)!}\sum_{(x_1,\dots ,x_{d-k})\subset\{u_1,\dots ,u_n\}}{\mathbb H}^k(\cap_{i=1}^{d-k}x_i)\right)\lambda^n (\dd(u_1,\dots ,u_n))$$ and the result follows.
\hfill $\Box $

\noindent We obtain asymptotic results in a special model.
Let \begin{equation}\label{ymod}Y=[0,b]^d\times\{2b\}\times\left\{e_i,\;i=1,\dots ,d\right\},\end{equation}
where $e_i$ are canonical unit vectors. A facet is a set $$y=(z,\phi )=\{s\in D; \max_i|z_i-s_i|\leq b\},$$ $z$ denotes the centre and $D$ is hyperplane with normal orientation $\phi .$ That means facets have the same fixed size and shape and any non-parallel facets intersect. In the case of $d-k$ facets with different orientations we have bounds
\begin{equation}\label{boh}b^k\leq {\mathbb H}^k(\cap_{i=1}^{d-k}x_i)\leq (2b)^k.\end{equation}
In the intensity $\lambda $ (\ref{inti}) we have $Q=\delta_{2b},$ the orientation distribution $V$ is uniform on $\{e_i,\;i=1,\dots ,d\}.$
Then for $u_i=(z_i,\phi_i),\,i=1,\dots ,n,$ it holds 
\begin{equation}\label{rdf}
\lambda^n(\dd(u_1,\dots, u_n))=\end{equation}$$\chi(z_1)\dd z_1\dots \chi(z_n)\dd z_n{1\over d^n}{\sum\dots\sum}_{n_1+\dots +n_d=n}\frac{n!}{n_1!\dots n_d!}
\otimes_{i=1}^d\delta_{e_i}^{n_i}(\dd\phi_1,\dots,\dd\phi_n).
$$ Denote $T=\int_{[0,b]^d}\chi (z)\dd z.$
\begin{remark}
In the special model it can be shown that for any $l \in \{2, \ldots ,d \}$ if $\nu_{i} = 0, i \neq l$ and $\nu_{l} \geq 0$, then $p \not\in L^{1}(P_{\eta})$, thus the conditions applied to parameters are not only sufficient, but also neccesarry conditions for density existence.
\end{remark}
\begin{lemma}\label{ll4}
Denote as in Lemma \ref{ll1}, formula (\ref{rods}) 
$$\rho_{d-k}(y_1,\dots ,y_{d-k};\mu_a^{(d-k)})=\frac{A(a)}{B(a)},\;k=0,\dots ,d-2.$$
Then we have for $y_1,\dots ,y_{d-k}$ with different orientations $$B(a)\geq \sum_{n=0}^\infty\left(\frac{aT}{d}\right)^n{\sum\dots\sum}_{n_1+\dots +n_d=n}\frac{1}{n_1!\dots n_d!}\times$$$$\times\exp\left(\nu_{d-k}(2b)^k\sum_{\{l_j\}_{j=1}^{d-k}\subset [d]}\prod_{i=1}^{d-k}n_{l_j}\right),$$

$$A(a)\leq \sum_{n=0}^\infty\left(\frac{aT}{d}\right)^n{\sum\dots\sum}_{n_1+\dots +n_d=n}\frac{1}{n_1!\dots n_d!}\times$$$$\times\exp\left(\nu_{d-k}b^k \left( \sum_{\{l_j\}_{j=1}^{d-k}\subset [d]}\prod_{i=1}^{d-k}n_{l_j}+\vartheta \right)\right),$$
where $\vartheta$ is the number of $(d-k)$-tuples (with different orientations) $$\{x_1,\dots ,x_{d-k}\}\subset\{u_1,\dots ,u_n,y_1,\dots ,y_{d-k}\}$$ such that there is at least one of $y_j$ among $x_1,\dots ,x_{d-k}.$\end{lemma}
\noindent{\bf Proof:} It follows from Lemma \ref{ll1}, inequality (\ref{boh}) and integration w.r.t. (\ref{rdf}). Here $n_j$ are interpreted as the number of $u_l,\;l=1,\dots n$ with orientation $e_j,\;j=1,\dots d.$ \hfill $\Box$
\begin{lemma}\label{ll5}Under the notation from Lemma \ref{ll4} it holds $$B(a)\geq \exp(\frac{1}{d}aT(d-k-1)),\;k=0,\dots ,d-2.$$ \end{lemma}
\noindent{\bf Proof:} Substituting $\beta=\frac{aT}{d}$ we have
\begin{gather}
\nonumber \sum_{n=0}^{\infty} \beta^{n} \sum_{n_{1}+ \ldots n_{d} = n} \frac{1}{n_{1}! \ldots n_{d}!} \exp(\nu_{d-k} (2b)^k \sum_{\{l_{j}\}_{j=1}^{d-k} \subset [d]} \prod_{j=1}^{d-k}n_{l_{j}}) =  \\
\nonumber \sum_{n_{1}=0}^{\infty} \ldots \sum_{n_{d}=0}^{\infty} \frac{\beta^{n_{1}+\ldots + n_{d}}}{n_{1}! \ldots n_{d}!}  \exp(\nu_{d-k} (2b)^k \sum_{\{l_{j}\}_{j=1}^{d-k} \subset [d]} \prod_{j=1}^{d-k}n_{l_{j}}) = \\
\nonumber \sum_{n_{1}=0}^{\infty} \ldots \sum_{n_{d-1}=0}^{\infty} \frac{\beta^{n_{1}+\ldots + n_{d-1}}}{n_{1}! \ldots n_{d-1}!}\exp\left( \beta\exp(\nu_{d-k} (2b)^k \sum_{\{l_{j}\}_{j=0}^{d-k-1} \subset [d-1]} \prod_{j=1}^{d-k-1} n_{l_{j}}) \right) \times  \\
\nonumber \times \exp\left( \nu_{d-k} (2b)^k \sum_{\{l_{j}\}_{j=1}^{d-k} \subset [d-1]} \prod_{i=1}^{d-k}n_{l_{j}}\right) \geq  \\
\nonumber \sum_{n_{1}=0}^{\infty} \ldots \sum_{n_{d-1}=0}^{\infty} \frac{\beta^{n_{1}+\ldots + n_{d-1}}}{n_{1}! \ldots n_{d-1}!} \exp( \nu_{d-k} (2b)^k \sum_{\{l_{j}\}_{j=1}^{d-k} \subset [d-1]} \prod_{j=1}^{d-k}n_{l_{j}}) \geq \\
\nonumber \sum_{n_{1}=0}^{\infty} \ldots \sum_{n_{d-k}=0}^{\infty} \frac{\beta^{n_{1}+\ldots + n_{d-k}}}{n_{1}! \ldots n_{d-k}!} \exp( \nu_{d-k} (2b)^k \sum_{\{l_{j}\}_{j=1}^{d-k} \subset [d-k]} \prod_{j=1}^{d-k}n_{l_{j}}) = \\
\nonumber \sum_{n_{1}=0}^{\infty} \ldots \sum_{n_{d-k}=0}^{\infty} \frac{\beta^{n_{1}+\ldots + n_{d-k}}}{n_{1}! \ldots n_{d-k}!} \exp( \nu_{d-k} (2b)^k  \prod_{i=1}^{d-k}n_{i}) =\\
\nonumber \sum_{n_{1}=0}^{\infty} \ldots \sum_{n_{d-k-1}=0}^{\infty}\frac{\beta^{n_{1}+\ldots + n_{d-k-1}}}{n_{1}! \ldots n_{d-k-1}!} \exp\left( \beta e^{\nu_{d-k} (2b)^{k}  \prod_{i=1}^{d-k-1}n_{i} } \right) \geq \\
\nonumber \sum_{n_{1}=0}^{\infty} \ldots \sum_{n_{d-k-1}=0}^{\infty}\frac{\beta^{n_{1}+\ldots + n_{d-k-1}}}{n_{1}! \ldots n_{d-k-1}!} = e^{\beta(d-k-1)} . 
\end{gather}\hfill $\Box $
\begin{lemma}\label{ll6}
There exists a constant $R<0$ independent of $a$ such that \begin{equation}\label{oddh} \frac{A(a)}{B(a)}\leq e^{Ra}\end{equation}\end{lemma}
\noindent{\bf Proof:} Substituting $\beta=\frac{aT}{d}$ we estimate $A(a)$  
\begin{equation}\label{prsu}
A(a)\leq\sum_{n=0}^{\infty}\beta^{n}\sum_{n_{1}+ \ldots +n_{d} = n} \frac{1}{n_{1}! \ldots n_{d}!} \exp ( \nu_{d-k} b^k C_{d-k,d}),
\end{equation}
where 
\begin{align}
\nonumber C_{l,m} &= \sum_{i=1}^{l} \sum_{\{j_{1}, \ldots , j_{i}\} \subset [m]} \prod_{k=1}^{i} n_{j_{k}}+1,\;1\leq l\leq m. \\  
\nonumber C_{l,m} &= C_{l-1,m-1}n_{m}+C_{l,m-1},\;2\leq l<m, \\
\nonumber C_{m,m} &= C_{m-1,m-1}(1+n_{m}).
\end{align}
Now we decrease the number of sums in (\ref{prsu}):
\begin{gather}
\nonumber \sum_{n_{1}=0}^{\infty} \ldots \sum_{n_{d}=0}^{\infty} \frac{\beta^{n_{1}+\ldots + n_{d}}}{n_{1}! \ldots n_{d}!} \exp ( \nu_{d-k} b^k C_{d-k,d}) = \\
\nonumber \sum_{n_{1}=0}^{\infty} \ldots \sum_{n_{d}=0}^{\infty} \frac{\beta^{n_{1}+\ldots + n_{d}}}{n_{1}! \ldots n_{d}!} \exp ( \nu_{d-k} b^k (C_{d-k-1,d-1}n_{d}+C_{d-k,d-1} ))= \\
\nonumber \sum_{n_{1}=0}^{\infty} \ldots \sum_{n_{d-1}=0}^{\infty} \frac{\beta^{n_{1}+\ldots + n_{d-1}}}{n_{1}! \ldots n_{d-1}!} \exp (\beta e^{\nu_{d-k} b^k C_{d-k-1,d-1}}  + \nu_{d-k} b^k C_{d-k,d-1} ) \leq \\
\nonumber \sum_{n_{1}=0}^{\infty} \ldots \sum_{n_{d-1}=0}^{\infty} \frac{\beta^{n_{1}+\ldots + n_{d-1}}}{n_{1}! \ldots n_{d-1}!} \exp ( \beta e^{\nu_{d-k} b^k p}  + \nu_{d-k} b^k C_{d-k,d-1} ) + \\ 
\nonumber \sum_{n_{1}=0}^{p} \ldots \sum_{n_{d-1}=0}^{p} \frac{\beta^{n_{1}+\ldots + n_{d-1}}}{n_{1}! \ldots n_{d-1}!} \exp (\beta e^{\nu_{d-k} b^k C_{d-k-1,d-1}}  + \nu_{d-k} b^k C_{d-k,d-1} ) - \\
\nonumber \sum_{n_{1}=0}^{p} \ldots \sum_{n_{d-1}=0}^{p} \frac{\beta^{n_{1}+\ldots + n_{d-1}}}{n_{1}! \ldots n_{d-1}!} \exp (\beta e^{\nu_{d-k} b^k p}  + \nu_{d-k} b^k C_{d-k,d-1} ). 
\end{gather}
In the last step $(p+1)^{d-1}$ terms of infinite series were changed. Then we estimate the exponential part of series from above by $\exp(\beta e^{\nu_{d-k}p})$.\\
Each term of both finite series tends to zero when divided by $e^{\beta(d-k-1)}$ (as $\beta$ tends to infinity), so that the whole series tends to zero. This step is repeated till we get
\begin{gather}
\nonumber \alpha_{p} \sum_{n_{1}=0}^{\infty} \ldots \sum_{n_{d-k}=0}^{\infty} \frac{\beta^{n_{1}+\ldots + n_{d-k}}}{n_{1}! \ldots n_{d-k}!} \exp ( \nu_{d-k} b^k C_{d-k,d-k} ) \leq \\
\nonumber \alpha_{p} \sum_{n_{1}=0}^{\infty} \ldots \sum_{n_{d-k}=0}^{\infty} \frac{\beta^{n_{1}+\ldots + n_{d-k}}}{n_{1}! \ldots n_{d-k}!} \exp \left( \nu_{d-k} b^k  \left[ \sum_{i=1}^{d-k-1}n_{i}n_{d-k}+ \sum_{i=1}^{d-k}n_{i} \right] \right) = \\
\nonumber \alpha_{p} \sum_{n_{1}=0}^{\infty} \ldots \sum_{n_{d-k-1}=0}^{\infty} \frac{(\beta e^{\nu_{d-k} b^k})^{n_{1}+\ldots + n_{d-k-1}}}{n_{1}! \ldots n_{d-k-1}!} \exp (\beta e^{\nu_{d-k} b^k (\sum_{i=1}^{d-k-1}n_{i}+1)}) \leq 
\end{gather}
\begin{gather}
\nonumber \alpha_{p} \sum_{n_{1}=0}^{\infty} \ldots \sum_{n_{d-k-1}=0}^{\infty} \frac{(\beta e^{\nu_{d-k} b^k})^{n_{1}+\ldots + n_{d-k-1}}}{n_{1}! \ldots n_{d-k-1}!} \exp (\beta e^{\nu_{d-k} b^k q})+ \\
\nonumber \alpha_{p} \sum_{n_{1}=0}^{q} \ldots \sum_{n_{d-k-1}=0}^{q} \frac{(\beta e^{\nu_{d-k} b^k})^{n_{1}+\ldots + n_{d-k-1}}}{n_{1}! \ldots n_{d-k-1}!} \exp (\beta e^{\nu_{d-k} b^k (\sum_{i=1}^{d-k-1}n_{i}+1)}) - \\
\nonumber \alpha_{p} \sum_{n_{1}=0}^{q} \ldots \sum_{n_{d-k-1}=0}^{q} \frac{(\beta e^{\nu_{d-k} b^k})^{n_{1}+\ldots + n_{d-k-1}}}{n_{1}! \ldots n_{d-k-1}!} \exp (\beta e^{\nu_{d-k} b^k q}),
\end{gather}
where $\alpha_{p} = \exp ( \beta e^{\nu_{d-k} b^k p}k ) $. Again we changed $(q+1)^{d-k-1}$ terms of infinite series and estimated exponential part of the series by $\exp(\beta e^{\nu_{d-k}q})$. Both finite series tend to zero when divided by $e^{\beta (d-k-1)}$ and the infinite series is equal to 
\begin{gather}
\nonumber \exp ( \beta e^{\nu_{d-k} b^k p}k +\beta e^{\nu_{d-k} b^k}(d-k-1) + \beta e^{\nu_{d-k} b^k q}).
\end{gather}
We have to select $p$ and $q$ so that
\begin{gather}
\nonumber R_1=e^{\nu_{d-k} b^k p}k +e^{\nu_{d-k} b^k}(d-k-1) + e^{\nu_{d-k} b^k q} - (d-k-1)<0,
\end{gather}
which is possible for any $\nu_{d-k}<0.$
Then we reverse the substitution $\beta=\frac{aT}{d}$ and the ratio $\frac{A(a)}{B(a)}$ tends to zero for $a\rightarrow\infty.$ Its convergence is not slower than $e^{Ra},$ where $R=\frac{TR_1}{d}.$ \hfill $\Box $
\begin{theorem}\label{tecko}For $a\rightarrow\infty$ it holds \begin{equation}\label{a}{\mathbb E} G_{d-k}(\mu^{(d-k)}_a)\rightarrow 0,\; k=0,\dots ,d-2.\end{equation}\end{theorem}
\noindent {\bf Proof:} In (\ref{znovu}) it suffices to show that for $y_1,\dots ,y_{d-k}$ all with distinct orientation we have \begin{equation}\label{roj}\lim_{a\rightarrow\infty}a^{d-k}\rho_{d-k}(y_1,\dots ,y_{d-k};\mu_a^{(d-k)})=0.\end{equation} This is the consequence of Lemmas \ref{ll5} and \ref{ll6}. Then by the Lebesgue dominance theorem formula (\ref{a}) follows. \hfill $\Box $
\begin{remark}
Since the functionals $G_{d-k}$ are nonnegative Theorem \ref{tecko} says that asymptotically the processes $\mu_a^{d-k}$ are degenerate (for any $\nu_{d-k}<0$) in the sense that there are no intersections of $(d-k)$-tuples, i.e. some orientations are missing. For $k=d-2$ e.g. all facets tend to be parallel (with any orientation). Additionally when the assertion of Theorem \ref{tecko} is valid for $\mu_a^{(d-k)}$ then it holds also for expectations ${\mathbb E} G_j(\mu^{(d-k)}_a) $ with $j=d-k+1,\dots ,d.$ \end{remark}
\begin{remark} The behavior of Metropolis-Hastings chain for simulation of a realization of $\mu_a^{d-k}$ with fixed parameter $a$ depends on $\nu_{d-k}$ in the model with finitely many orientations. For $|\nu_{d-k}|$ large typically it converges quickly to a realization with some missing orientations. However when $\nu_{d-k}$ is close to zero it does not converge quickly at all. This property does not take place when the orientation distribution is absolutely continuous w.r.t. spherical Lebesgue measure.\end{remark}

\section{$U$-statistics of the order smaller than the submodel}
In this section we continue to study the model from previous subsection, see (\ref{ymod}). Asymptotic moments (when $a\rightarrow\infty $) of functionals
$G_{d-k}(\x)$ in the submodel $\mu_{a}^{(d)},$ where $k=1,\ldots,d-1;\; d \geq 3,$ will be investigated. Here the vanishing property (\ref{a}) is not expected, see the following Table where crosses mean that the expected value is non-zero. 

\medskip

\begin{center}
 \renewcommand{\arraystretch}{1.5}
\begin{tabular}{|c|c|c|c|c|c|}
\hline $a\rightarrow\infty$ & \multicolumn{5}{|c|}{\textit{U}-statistics $G_j$} \\
\hline Submodel & $\mathbb{E}G_{d}$ & $\mathbb{E}G_{d-1}$ & \ldots  & $\mathbb{E} G_{2}$  & $\mathbb{E}G_{1}$\\
\hline $\mu^{(2)}_{a}$ & $0$ &  $0$  & \ldots & $0$  & $\times$ \\
\hline $\mu^{(3)}_{a}$ & $0$ &  $0$  & \ldots & $\times$  & $\times$ \\
\hline \vdots & \vdots  &   \vdots & & \vdots &\vdots \\

\hline 
$\mu^{(d-1)}_{a}$ &  $0$  & $0$ &\ldots & $\times$ & $\times$ \\
\hline
$\mu^{(d)}_{a}$ &  $0$  &  $\times$  & \ldots & $\times$  & $\times$\\
\hline
\end{tabular}
\end{center}

\medskip

As in the previous Section we will need first some limits of correlation functions.
\subsection{The limit of correlation functions}
\begin{proposition}\label{teh50}a) For any $y_1,\dots ,y_{d-k}\in Y$ with different orientations we have \begin{equation}\label{dals}
\lim_{a\rightarrow\infty}\rho_{d-k}(y_{1},\ldots,y_{d-k};\mu_{a}^{(d)}) = \frac{k}{d},\;k=1,\ldots,d-1.\end{equation}
b) For all sets of $2(d-k)$ arguments which fullfil that each facet of $\{ y_{1},\ldots,y_{d-k} \}$ has different orientation and the same applies to $\{ y_{d-k+1},\ldots,y_{2(d-k)} \}$ and there is $l$ common orientations within these two groups, then \begin {equation}\label{teh52}\lim_{a\rightarrow\infty}\rho_{2(d-k)}(y_{1},\ldots,y_{2(d-k)};\mu_{a}^{(d)})=\frac{2k-d+l}{d}.\end{equation}
c) Consider $2(d-k)-1$ arguments such that each facet of $\{ y_{1},\ldots,y_{d-k} \}$ has different orientation and the same applies to $\{ y_{d-k+1},\ldots,y_{2(d-k)-1},y_1 \}.$ Let us have $l$ common orientations $e_{2},\ldots,e_{l+1}$ within these two groups of arguments not including orientation $e_{1}$ of $y_{1}$, then
\begin {equation}\label{teh53}\lim_{a\rightarrow\infty}\rho_{2(d-k)-1}(y_{1},\ldots,y_{2(d-k)-1};\mu_{a}^{(d)})=\frac{2k-d+l-1}{d}.\end{equation}
\end{proposition}
\noindent{\bf Proof} 
a) Consider fixed facets $y_{1}, \ldots, y_{d-k}$ each associated with a different normal orientation vector, without loss of generality $e_{1}, \ldots, e_{d-k}.$ Then the correlation function can be expressed analogously to Lemma \ref{ll1} as follows:
$$
\rho_{d-k}(y_{1},\ldots,y_{d-k};\mu_{a}^{(d)}) = \frac{\tilde{A}(a)}{\tilde{B}(a)},\qquad\tilde{A}(a)=$$
$$= \sum_{n=0}^{\infty}\frac{a^n}{n!}\int_{Y} \ldots \int_{Y} exp \left( \frac{\nu_{d}}{d!} \sum_{\substack{(x_{1},\ldots,x_{d}) \subset\\ \{u_{1},\ldots,u_{n},y_{1},\ldots,y_{d-k} \} }} \mathbb{H}^{0}(\cap_{i=1}^{d}x_{i}) \right)   \lambda^n(\dd(u_{1}, \ldots , u_{n})),$$
$\tilde{B}(a)=$$$= \sum_{n=0}^{\infty}\frac{a^n}{n!}\int_{Y} \ldots \int_{Y} exp \left( \frac{\nu_{d}}{d!} \sum_{(x_{1},\ldots,x_{d}) \subset \{u_{1},\ldots,u_{n} \} } \mathbb{H}^{0}(\cap_{i=1}^{d}x_{i}) \right)\lambda^n(\dd(u_{1}, \ldots , u_{n})).
$$
It holds $\tilde{A}(a)/\tilde{B}(a)=A(a)/B(a),$ where using (\ref{rdf}) we write
$A(a)= e^{-\frac{aT}{d}(d-1)}\times $ 
\begin{equation}\label{ab}\times \sum_{n=0}^{\infty} \left( \frac{aT}{d} \right)^n \sum_{n_{1}+ \ldots + n_{d}=n} \frac{1}{n_{1}! \ldots n_{d}!} \exp \left(   \nu_{d} \sum_{D\subset [d-k]} \prod_{l \in [d]\setminus D} n_{l} \right), \end{equation}
$$B(a)=
\nonumber e^{-\frac{aT}{d}(d-1)} \sum_{n=0}^{\infty} \left( \frac{aT}{d} \right)^n \sum_{n_{1}+ \ldots + n_{d}=n} \frac{1}{n_{1}! \ldots n_{d}!} \exp \left(   \nu_{d}\prod_{l = 1}^{d} n_{l} \right). $$
From Lemmas \ref{lemB2} and \ref{lemA2} below, (\ref{dals}) follows. 

\noindent b) We need to compute $\rho_{2(d-k)}(y_{1},\ldots,y_{2(d-k)};\mu_{a}^{(d)})$ for all sets of arguments which fullfil that each facet of $\{ y_{1},\ldots,y_{d-k} \}$ has different orientation and the same applies to $\{ y_{d-k+1},\ldots,y_{2(d-k)} \}$. Without loss of generality  we need to consider only situations where facets $ y_{1},\ldots,y_{d-k}$  have orientations $e_{1},\ldots,e_{d-k}$ and facets $y_{d-k+1},\ldots,y_{2(d-k)}$ have orientations $e_{1},\ldots,e_{l},e_{d-k+1},\ldots e_{2d-2k-l}$, where $l =\max(d-2k,0),\ldots,d-k $, probability of both sets of facets having $l$ common orientations can be expressed by multinomial coefficient
\begin{equation}\label{ddbb}
\frac{(d-k)!^{2}}{d^{2(d-k)}} \binom{d}{l,d-k-l,d-k-l,2k-d+l}.
\end{equation}
It holds
\begin{gather}
\nonumber \lim_{a\rightarrow\infty}\rho_{2(d-k)}(y_{1},\ldots,y_{2(d-k)};\mu_{a}^{(d)}) = \lim_{a\rightarrow\infty}\frac{e^{-a(d-1)}}{d} \sum_{n=0}^{\infty} \frac{a^{n}}{n!}\times  \\
\nonumber \times	\int_{Y^n} \exp \left( \frac{\nu_{d}}{d!} \sum_{ (x_{1},\ldots,x_{d} ) \subset \{u_{1},\ldots,u_{n},y_{1},\ldots,y_{2(d-k)} \} } \mathbb{H}^{0}(\cap_{i=1}^{d}x_{i}) \right)  \lambda^n(\dd(u_{1}, \ldots,u_{n})),
\end{gather}
where we get $\frac{e^{-a(d-1)}}{d}$ from calculations of the limit $B(a)$ in Lemma \ref{lemB2}. For selected $l$ the limit of the correlation function can be further expressed in the form
\begin{gather}
\label{lt}\lim_{a\rightarrow\infty}\frac{e^{-a(d-1)}}{d} \sum_{n_{1}=0}^{\infty} \ldots \sum_{n_{d}=0}^{\infty}
\frac{a^{n_{1}+ \ldots + n_{d}}}{n_{1}! \ldots n_{d}!} \\
\nonumber \times \exp \left( \nu_{d} \sum_{D \subset [2d-2k-l]} \prod_{m \in [d]\setminus D} n_{m} 2^{|D \cap \{1,\ldots, l \}|} \right).
\end{gather}
The expression in the exponent can be bounded from both sides
\begin{gather}
\nonumber \nu_{d} 2^{l} \sum_{D \subset [2d-2k-l]} \prod_{m \in [d]\setminus D} n_{m} \leq \\
\nonumber \nu_{d} \sum_{D \subset [2d-2k-l]} \prod_{m \in [d]\setminus D}n_{m} 2^{|D \cap \{1,\ldots, l \}|} \leq \\
\nonumber \nu_{d} \sum_{D \subset [2d-2k-l]} \prod_{m \in [d]\setminus D}n_{m}.
\end{gather}
For the both bounding series the limit (\ref{lt}) is equal to $\frac{2k-d+l}{d}$ since the series are in the form as $A(a)$ in Lemma \ref{lemA2}.

\noindent c) In the previous case it can be seen the numerator is equal to number of unused orientations among facets $y_{1},\ldots,y_{2(d-k)}$ and the same applies to this calculation except that there is one extra common orientation of $x_{1}$. \hfill $\Box $

\begin{remark}\label{rmm}
In b) the limit is equal to zero for $l=d-2k.$ In this case all orientations up to the order $d$ of the submodel are exhausted. This corresponds to the situation from Theorem \ref{tecko} where in (\ref{roj}) also all orientations up to the order $d-k$ of the submodel are exhausted by $y_1,\dots ,y_{d-k}$ and the correlation function tends to zero. In the opposite case (in Theorem \ref{thvar} when $l>d-2k$ is admissible) the limit of correlation function is nonzero.
\end{remark} 
In the following two lemmas we write for simplicity $\nu $ instead of $\nu_d.$
\begin{lemma}\label{lemB2} It holds $\lim_{a\rightarrow\infty}B(a)=d.$\end{lemma}
\noindent{\bf Proof:} We want to examine series in form
\begin{equation}\label{wes}
\sum_{n_{1}=0}^{\infty} \ldots \sum_{n_{d}=0}^{\infty}
\frac{a^{n_{1}+ \ldots + n_{d}}}{n_{1}! \ldots n_{d}!} \exp ( \nu n_{1} \ldots n_{d} -a(d-1))= \end{equation}
$$\sum_{n_{1}=0}^{\infty} \ldots \sum_{n_{d-1}=0}^{\infty}
\frac{a^{n_{1}+ \ldots + n_{d-1}}}{n_{1}! \ldots n_{d-1}!} \exp ( ae^{ \nu n_{1} \ldots n_{d-1}} -a(d-1)),
$$
where $\nu<0$ and $d\geq 3$. This form fits the $B(a)$ in (\ref{ab}) after substituting $a$ for $\frac{aT}{d}.$ 
We divide indices in the sums into three subsets
\begin{equation}\label{rozdel}
D_1 = \{ n_{1} \geq \sqrt{a} \wedge  \ldots \wedge n_{d-1} \geq  \sqrt{a} \},\end{equation}
$$D_2 = \{ n_{1} = 0 \vee  \ldots \vee n_{d-1} = 0\},$$
$$D_3 = \{ n_{1} < \sqrt{a} \vee  \ldots \vee n_{d-1} <  \sqrt{a} \} \setminus D_2.$$
Then we use one of Chernoff's bounds for Poisson distribution, which says that it holds
\begin{align}
\nonumber \sum_{k=0}^{t} \frac{s^{k}}{ k!} e^{-s} \leq e^{-s} \frac{ (es)^t}{ t^t}, t<s.
\end{align}
Firstly, we sum over $D_2$. Let $n_{1}=0,$ then from (\ref{wes}) we have
\begin{align}
\nonumber \sum_{n_{2}=0}^{\infty} \ldots \sum_{n_{d-1}=0}^{\infty}
\frac{a^{n_{2}+ \ldots + n_{d-1}}}{n_{2}! \ldots n_{d-1}!} \exp ( a -a(d-1)) = e^{ a -a(d-1)+a(d-2)}=1.
\end{align}
We continue by induction with $n_1\neq 0,\dots ,n_{k-1}\neq 0,\,n_k=0$ where we get 
$$
\sum_{n_{1}=1}^{\infty} \ldots \sum_{n_{k-1}=1}^{\infty} \sum_{n_{k+1}=0}^{\infty} \ldots \sum_{n_{d-1}=0}^{\infty}
\frac{a^{n_{1}+ \ldots + n_{k-1} + n_{k+1} + \ldots + n_{d-1}}}{n_{1}!  \ldots n_{k-1}!  n_{k+1}! \ldots n_{d-1}!} \exp ( a -a(d-1)) =$$$$
\nonumber =e^{-a(k-1)}(e^{a}-1)^{k-1} \rightarrow 1.
$$
By calculating all $d-1$ options we explored all combinations of indices in $D_2$ and we conclude that the sum over $D_2$ tends to $d-1$.

Secondly, we sum over $D_3$.
\begin{gather}
\nonumber  \sum_{ \{ n_{1},\ldots,n_{d-1} \} \in D_3}
\frac{a^{n_{1}+ \ldots + n_{d-1}}}{n_{1}! \ldots n_{d-1}!} \exp ( ae^{ \nu n_{1} \ldots n_{d-1}} -a(d-1)) \leq \\
\nonumber \sum_{ \{ n_{1},\ldots,n_{d-1} \} \in D_3}
\frac{a^{n_{1}+ \ldots + n_{d-1}}}{n_{1}! \ldots n_{d-1}!} \exp (ae^{ \nu } -a(d-1)) \leq \\
\nonumber e^{ae^{ \nu }} \sum_{ \{n_{1},\ldots,n_{d-1} \} \in D_2  \cup D_3} \frac{a^{n_{1}+ \ldots + n_{d-1}}}{n_{1}! \ldots n_{d-1}!} \exp ( -a(d-1)) \leq \\
\nonumber e^{ae^{ \nu }} \left(  \sum_{n_{1}=0}^{ \lfloor \sqrt{a} \rfloor }\frac{a^{n_{1}}}{n_{1}!}e^{-a} + \ldots + \sum_{n_{d-1}=0}^{\lfloor \sqrt{a} \rfloor
}\frac{a^{n_{d-1}}}{n_{d-1}!}e^{-a} \right) \leq \\
\nonumber e^{ae^{ \nu }}(d-1) e^{-a} \frac{ (ea)^{\lfloor \sqrt{a} \rfloor }}{(\lfloor \sqrt{a} \rfloor)^{\lfloor \sqrt{a} \rfloor}},
\end{gather}
where we used the principle of inclusion and exclusion (for probabilities of Poisson distribution) and the Chernoff's bound. Then we examine logarithm of the previous term
\begin{align}
\nonumber a(e^{ \nu }-1)  + \log(d-1)+ \lfloor \sqrt{a} \rfloor \log(ea)-\lfloor \sqrt{a} \rfloor \log(\lfloor \sqrt{a} \rfloor).
\end{align}
which tends to $- \infty$ and so the sum over $D_3$ tends to $0$, when $a\rightarrow\infty $.\\
Finally, we examine the sum over indices in $D_1$ in order to show that it tends to $1.$ Let $\varepsilon>0$ be arbitrarily chosen. We choose $\gamma_1$, so that $$|\exp(ae^{\nu ( \lceil \sqrt{a} \rceil )^{d-1}}) -1 | < \varepsilon, \forall a \geq \gamma_1.$$ In the next step we choose $\gamma_2$, which fullfills the following condition:
\begin{align}
\nonumber e^{-a} \frac{ (ea)^{\lfloor \sqrt{a} \rfloor }}{(\lfloor \sqrt{a} \rfloor)^{\lfloor \sqrt{a} \rfloor}} < \varepsilon,\forall a \geq \gamma_2.
\end{align}
Then we estimate the series from both sides 
\begin{equation}\label{polic}
\sum_{n_{1}=  \lceil \sqrt{a}  \rceil }^{\infty}  \ldots \sum_{n_{d-1}=  \lceil \sqrt{a}  \rceil }^{\infty}\frac{a^{n_{1}+ \ldots + n_{d-1}}}{n_{1}! \ldots n_{d-1}!} e^{-a(d-1)} \leq  \end{equation}\begin{gather} 
\nonumber \sum_{n_{1}=  \lceil \sqrt{a}  \rceil }^{\infty}  \ldots \sum_{n_{d-1}=  \lceil \sqrt{a}  \rceil }^{\infty}\frac{a^{n_{1}+ \ldots + n_{d-1}}}{n_{1}! \ldots n_{d-1}!} e^{ae^{\nu n_{1} \ldots n_{d-1}}-a(d-1)} \leq  \\
\nonumber \sum_{n_{1}=  \lceil \sqrt{a}  \rceil }^{\infty}  \ldots \sum_{n_{d-1}=  \lceil \sqrt{a}  \rceil }^{\infty}\frac{a^{n_{1}+ \ldots + n_{d-1}}}{n_{1}! \ldots n_{d-1}!} e^{-a(d-1)} (1+\varepsilon) \leq (1+\varepsilon). 
\end{gather}
Investigating the lower bound it can be seen that 
\begin{gather}
\nonumber \sum_{n_{1}=  \lceil \sqrt{a}  \rceil }^{\infty}  \ldots \sum_{n_{d-1}=  \lceil \sqrt{a}  \rceil }^{\infty}\frac{a^{n_{1}+ \ldots + n_{d-1}}}{n_{1}! \ldots n_{d-1}!} e^{-a(d-1)} = \\
\nonumber 1-\sum_{ \{ n_{1}, \ldots n_{d-1} \} \in D_2 \cup D_3}\frac{a^{n_{1}+ \ldots + n_{d-1}}}{n_{1}! \ldots n_{d-1}!} e^{-a(d-1)} \geq \\
\nonumber 1-(d-1) \sum_{n=0}^{ \lfloor \sqrt{a} \rfloor}\frac{a^{n}}{n!}e^{-a} \geq 1 - (d-1)e^{-a} \frac{ (ea)^{\lfloor \sqrt{a} \rfloor }}{(\lfloor \sqrt{a} \rfloor)^{\lfloor \sqrt{a} \rfloor}} \geq 1 - (d-1)\varepsilon,
\end{gather}
where we used again the principle of inclusion and exclusion and the Chernoff's bound. Thus  the sum over $D_1$ can be enclosed by bounds which are arbitrarily close to 1.
We conclude that the overall sum (\ref{wes}) tends to $d$.\hfill $\Box $

\begin{lemma}\label{lemA2}It holds $\lim_{a\rightarrow\infty}A(a)=k.$\end{lemma}
\noindent{\bf Proof:}
We examine the expression for $A(a)$ in (\ref{ab}) substituting $\frac{aT}{d}=a,$ in the form 
\begin{equation}\label{was}\sum_{n_{1}=0}^{\infty} \ldots \sum_{n_{d}=0}^{\infty}
\frac{a^{n_{1}+ \ldots + n_{d}}}{n_{1}! \ldots n_{d}!} \exp ( \nu Q_{d-k,d} -a(d-1)) ,\; 1 \leq k\leq d-1 ,\end{equation} where 
$$Q_{s,t}= \sum_{F\subset [s]} \prod_{l \in [t]\setminus F} n_{l}  ,\;  0 \leq s \leq t,$$$
\prod_{l \in \emptyset} n_{l} = 1.$
 It holds
\begin{align}
\nonumber Q_{s,t} = Q_{s,t-1}n_{t}, s<t, \\
\nonumber Q_{s,t} = Q_{s,s}n_{s+1} \ldots n_{t}, s<t,  \\
\nonumber Q_{t,t} = Q_{t-1,t-1}(n_{t}+1), 
\end{align}
and it follows
$$\sum_{n_{1}=0}^{\infty} \ldots \sum_{n_{d}=0}^{\infty}
\frac{a^{n_{1}+ \ldots + n_{d}}}{n_{1}! \ldots n_{d}!} \exp ( \nu Q_{d-k,d} -a(d-1)) = 
$$
\begin{equation}\label{csc} =\sum_{n_{1}=0}^{\infty} \ldots \sum_{n_{d-1}=0}^{\infty}
\frac{a^{n_{1}+ \ldots + n_{d-1}}}{n_{1}! \ldots n_{d-1}!} \exp ( ae^{\nu Q_{d-k,d-1}} -a(d-1)).
\end{equation}
For the same subsets of indices as in (\ref{rozdel}) we use the fact that current series are bounded from above by the corresponding ones.
Therefore the arguments from Lemma \ref{lemB2} remain the same for $D_3,$ and also for $D_1$  
since (\ref{csc}) is ordered between the first and second expression in (\ref{polic}).

The subset $D_2$ has to be considered where we proceed analogously to the proof of Lemma 7, but in reverse order from $n_{d-1}$ to $n_{1}.$ The fact that in all terms of $Q_{d-k,d-1}$ there are factors $n_{d-1},\dots n_{d-k+1}$ is used. Setting $n_{d-t}$ to zero we get
\begin{gather}
\nonumber \sum_{n_{1}=0}^{\infty} \ldots \sum_{n_{d-t-1}=0}^{\infty}\sum_{n_{d-t+1}=1}^{\infty} \ldots \sum_{n_{d-1}=1}^{\infty}
\frac{a^{n_{1}+ \ldots + n_{d-t-1} + n_{d-t+1}+ \ldots  + n_{d-1}}}{n_{1}! \ldots n_{d-t-1}!  n_{d-t+1}!\ldots n_{d-1}!} \times \\ \times \exp ( a -a(d-1)) = 
\nonumber  (e^a-1)^{t-1}e^{a(d-t-1)}e^{a(d-2)}\rightarrow 1,\; k > t, \\
\nonumber \sum_{n_{1}=0}^{\infty} \ldots \sum_{n_{d-t-1}=0}^{\infty}\sum_{n_{d-t+1}=1}^{\infty} \ldots \sum_{n_{d-1}=1}^{\infty}
\frac{a^{n_{1}+ \ldots + n_{d-t-1} + n_{d-t+1}+ \ldots  + n_{d-1}}}{n_{1}! \ldots n_{d-t-1}!  n_{d-t+1}!\ldots n_{d-1}!} \times \\
\times  \exp ( ae^{\nu} -a(d-1)) = 
\nonumber (e^a-1)^{t-1}e^{a(d-t-1)}e^{a\nu -a(d-1)} \rightarrow 0,\; k \leq t,
\end{gather}
where the second expression is an upper bound. Thus the sum over $D_2$ tends to $k-1$ and we conclude that the overall sum (\ref{was}) tends to $k.$\hfill $\Box $

\subsection{The asymptotics of moments}
In this subsection first and second moments of functionals
$G_{d-k}(\x)$ in the submodel $\mu_{a}^{(d)},$ where $k=1,\ldots,d-1,\; d \geq 3,$ are studied. We will need the following Lemma.
\begin{lemma}\label{vcr}
Let two sums of finite length $\sum_{i=1}^{K} n^{(a)}_{i}$ and $\sum_{i=1}^{L} m^{(a)}_{i}$
fullfill
$$\lim_{a \rightarrow \infty} \sum_{i=1}^{K} n^{(a)}_{i}=\lim_{a \rightarrow \infty}\sum_{i=1}^{L} m^{(a)}_{i},$$
$n^{(a)}_{i}= \frac{c^{(a)}_{i,1}}{c^{(a)}_{i,2}}$, $m^{(a)}_{i}= \frac{d^{(a)}_{i,1}}{d^{(a)}_{i,2}},$
$$\lim_{a \rightarrow \infty} c^{(a)}_{i,1} = c_{i,1} \in (0,\infty),\;\lim_{a \rightarrow \infty} c^{(a)}_{i,2} = c_{i,2} \in (0,\infty),$$
$$\lim_{a \rightarrow \infty} d^{(a)}_{i,1} = d_{i,1} \in (0,\infty),\;\lim_{a \rightarrow \infty} d^{(a)}_{i,2} = d_{i,2} \in (0,\infty).$$
Let there exist $a_{0}>0,\;\gamma <0$ and $\zeta > 0$ such that for all $a\geq a_{0}$ and all $i$ it holds $$| c^{(a)}_{i,1}- c_{i,1}  | < \zeta e^{\gamma a},\; |c^{(a)}_{i,2}- c_{i,2}  | <\zeta e^{\gamma a},$$
$$| d^{(a)}_{i,1}- d_{i,1}  | <\zeta e^{\gamma a},\; | d^{(a)}_{i,2}- d_{i,2}  | <\zeta e^{\gamma a},$$
then for $p>0$ it holds
$$\lim_{a \rightarrow \infty} a^{p} \left\vert \sum_{i=1}^{K} \frac{c^{(a)}_{i,1}}{c^{(a)}_{i,2}} - \sum_{i=1}^{L} \frac{d^{(a)}_{i,1}}{d^{(a)}_{i,2}} \right\vert = 0.$$
\end{lemma}
{\bf Proof:} It holds
$$\lim_{a \rightarrow \infty} a^{p} \left\vert \sum_{i=1}^{K} \frac{c^{(a)}_{i,1}}{c^{(a)}_{i,2}} - \sum_{i=1}^{L} \frac{d^{(a)}_{i,1}}{d^{(a)}_{i,2}} \right\vert \leq 
 \lim_{a \rightarrow \infty} a^{p} \left\vert \sum_{i=1}^{K} \frac{c^{(a)}_{i,1}}{c_{i,2}-e^{\gamma a}} - \sum_{i=1}^{L} \frac{d^{(a)}_{i,1}}{d_{i,2}+e^{\gamma a}} \right\vert = $$$$
=  \left\vert \sum_{i=1}^{K} \frac{\lim_{a \rightarrow \infty} a^{p} c^{(a)}_{i,1}}{c_{i,2}} - \sum_{i=1}^{L} \frac{ \lim_{a \rightarrow \infty}  a^{p} d^{(a)}_{i,1}}{d_{i,2}} \right\vert \leq $$$$
\left\vert \sum_{i=1}^{K} \frac{\lim_{a \rightarrow \infty} a^{p} \zeta e^{\gamma a} }{c_{i,2}}\right\vert  +  \left\vert \sum_{i=1}^{L} \frac{ \lim_{a \rightarrow \infty}  a^{p} \zeta e^{\gamma a}}{d_{i,2}} \right\vert + \lim_{a \rightarrow \infty}a^{p} \left\vert \sum_{i=1}^{K} \frac{c_{i,1}}{c_{i,2}}  -\sum_{i=1}^{L} \frac{d_{i,1}}{d_{i,2}}  \right\vert =0.$$
\hfill $\Box $

We denote 
\begin{equation}\label{ik}I_k=\int_{([0,b]^{d})^{d-k}} \mathbb{H}^{k}(\cap_{i=1}^{d-k}(s_{i},2b,e_{i})) \chi(s_{1})\ldots \chi(s_{d-k})\dd s_{1}\ldots\dd s_{d-k},\end{equation}
\begin{align}
\nonumber  I^{\prime}_{k} = \int_{([0,b]^{d})^{2(d-k)-1}} {\mathbb H}^{k} (\cap_{i=1}^{d-k} (s_{i},2b,e_{i})) {\mathbb H}^{k} (\cap_{i=2}^{d-k} (s_{i+d-k-1},2b,e_{i}) \cap (s_{1},2b,e_{1}) ) \times \\
\nonumber \times \chi(s_{1}) \dd s_{1},\ldots,\chi(s_{2(d-k)-1}) \dd s_{2(d-k)-1},
\end{align}
 for facets $y_i=(s_{i},2b,e_{i}),\,i=1,\dots ,d-k,$ with different orientations. Because of symmetry $I_k,\;I^{\prime}_{k} $ do not depend on the choice of these orientations. From (\ref{boh}) it is \begin{equation}\label{nerz}0<T^{d-k}b^k\leq I_k\leq T^{d-k}(2b)^k.\end{equation}
\begin{theorem}\label{thvar}It holds for $k=1,\dots ,d-1:$
\begin{equation}
\label{teh5} \lim_{a\rightarrow\infty} \frac{G_{d-k}(\mu^{(d)}_a) }{a^{d-k}} = \frac{I_k }{d^{d-k}} \binom{d-1}{d-k} ~ a.s.,
\end{equation} 
\begin{equation}\label{teh6} \lim_{a \rightarrow \infty} \frac{1}{a^{2(d-k)-1}} \var G_{d-k}(\mu_{a}^{(d)}) = I^{\prime}_{k} \frac{(d-1)(d-k)^{2}}{d^{2(d-k)-1}} \binom{d-2}{d-k-1}^{2}. \end{equation}
\end{theorem}
\noindent{\bf Proof:}
First we show that \begin{equation}\label{dalsi}
\lim_{a\rightarrow\infty}  \frac{{\mathbb E} G_{d-k}(\mu^{(d)}_a)}{a^{d-k}} = \frac{I_k }{d^{d-k}} \binom{d-1}{d-k}.\end{equation}
From (\ref{effa}) and (\ref{hrm}) we have \begin{equation}\label{lote}{\mathbb E}{G_{d-k}(\mu_{a}^{(d)})} =
\end{equation}
 $$=\frac{a^{d-k}}{(d-k)!} \int_{Y^{d-k}} \mathbb{H}^{k}(\cap_{i=1}^{d-k}y_{i}) \rho_{d-k}(y_{1},\ldots,y_{d-k};\mu_{a}^{(d)}) \lambda^{d-k}(\dd (y_{1}, \ldots ,y_{d-k})).$$
Finally in order to obtain (\ref{dalsi}) we use (\ref{dals}), the fact that $\frac{d(d-1)...(d-k+1)}{d^{d-k}}$ is the probability of facets having different orientations and the Lebesgue dominance theorem in (\ref{lote}).

Further we evaluate the second moment of the $U$-statistics according to Theorem 1 using notation $$\bar{\mathbb H}^{k}(y_1,\dots ,y_{d-k})={\mathbb H}^{k} (\cap_{i=1}^{d-k}y_i).$$ By $\simeq$ we stress that only terms with higher or equal rate of convergence than $a^{2(d-k)-1}$ are expressed:
 \begin{equation}\label{dvac} 
\E G^{2}_{d-k}(\mu_{a}^{(d)}) =\sum_{\sigma \in \prod_{d-k,d-k}} \frac{a^{ |\sigma| }}{ ((d-k)!)^{2}}\times\end{equation}$$ \int_{Y^{|\sigma|}} \left(\bigotimes_{i=1}^{2} \bar{\mathbb H}^{k}  \right)_{\sigma}(y_{1}, \ldots,y_{|\sigma|})\rho_{|\sigma |}(y_1,\dots ,y_{|\sigma |};\mu_a^{(d)}) \lambda (\dd (y_{1},\ldots,y_{| \sigma|}))$$$$
\simeq J_1+J_2,$$ where $$J_1=\frac{a^{2(d-k)}}{((d-k)!)^{2}} \int_{Y^{2(d-k)}} {\mathbb H}^{k} (\cap_{i=1}^{d-k} y_{i} )  {\mathbb H}^{k} (\cap_{i=d-k+1}^{2(d-k)} y_{i} )\times $$$$\times\rho_{2(d-k)}(y_1,\dots ,y_{2(d-k)};\mu_a^{(d)})  \lambda^{2(d-k)} (\dd (y_{1},\ldots,y_{2(d-k)})),$$ 
$$J_2=\frac{a^{2(d-k)-1} (d-k)^{2}}{((d-k)!)^{2}}\int_{Y^{2(d-k)-1}} {\mathbb H}^{k} (\cap_{i=1}^{d-k} y_{i} )  {\mathbb H}^{k} (\cap_{i=d-k+1}^{2(d-k)-1} y_{i} \cap y_{1} )\times $$$$\times \rho_{2(d-k)-1}(y_1,\dots ,y_{2(d-k)-1};\mu_a^{(d)})  \lambda^{2(d-k)} (\dd (y_{1},\ldots,y_{2(d-k)-1})).$$
Here $(d-k)^{2}$ is the number of possible selections of a common element in the both functions. The first term of (\ref{dvac}) tends for $a\rightarrow\infty $ to $$\frac{I^{2}_{k}}{d^{2(d-k)}} \binom {d-1} {d-k} ^{2}.$$ From the Lebesgue dominance theorem in $J_1,$ using (\ref{rdf}), (\ref{ddbb}), (\ref{teh52}):
\begin{gather}
\nonumber \lim_{a\rightarrow\infty}\frac{{\mathbb E} G^{2}_{d-k}(\mu^{(d)}_a)}{a^{2(d-k)}} =\frac{I_k^{2}}{d^{2(d-k)}} \sum_{l = 0}^{d-k} \binom{d}{l,d-k-l,d-k-l,2k-d+l}\frac{2k-d+l}{d} \\
 \nonumber = \frac{I_k^{2}}{d^{2(d-k)}} \sum_{l = 0}^{d-k} \binom{d-1}{l,d-k-l,d-k-l,2k-d+l-1}= 
\frac{I_k^{2}}{d^{2(d-k)}}  \binom{d-1}{d-k}^{2}. 
\end{gather}
Here the last equation holds because both sides represent (up to a multiplicative constant) number of possibilities how to choose two sets of distinct $d-k$ elements from total $d-1$ elements.
This term cancels out with the squared expectation limit in (\ref{dalsi}), therefore $$\lim_{a\rightarrow\infty}\frac{\var G_{d-k}(\mu^{(d)}_a)}{a^{2(d-k)}}=0$$ and together with (\ref{dalsi}) the assertion (\ref{teh5}) of the theorem follows.

Because of this result we try further standardization by $a^{d-k-\frac{1}{2}},$ it is 
$$\frac{1}{a^{2(d-k)-1}} \var G_{d-k}(\mu_{a}^{(d)})\simeq\frac{1}{a^{2(d-k)-1}} (J_1+J_2-(\E G_{d-k}(\mu_{a}^{(d)}))^2).$$
First we shall prove that \begin{equation}\label{lmmm}\lim_{a\rightarrow\infty}a\left(\frac{J_1-(\E G_{d-k}(\mu_{a}^{(d)}))^2}{a^{2(d-k)}}\right)=0,\end{equation} which follows from Lemma \ref{vcr}. Here the two sums in the Lemma \ref{vcr} represent the two terms in the numerator of (\ref{lmmm}), which in our case can be expressed in such a way that each summand is a correlation function of some orientation configuration multiplied by $I_{k}$ and by the probability of the occurence of such configuration. The convergence rate properties follow from the calculations of the denominator of correlation function in Lemma \ref{lemB2}, where all parts converge at least at exponential rate to some positive value (we omit all zero values). The same holds for numerator of correlation function which is calculated in the same manner in Lemma \ref{lemA2}.

In the term $J_2$ of (\ref{dvac}) we need to consider number of common orientations between both functions. Let us have $l$ common orientations $e_{2},\ldots,e_{l+1}$ not including orientation $e_{1}$ of $x_{1}$, probability of the configuration with this number of common orientations is 
\begin{align}
\nonumber \frac{d}{d^{2(d-k)-1}} \binom{d-1}{l,d-k-l-1,d-k-l-1,2k+l-d+1}
\end{align}
and the corresponding correlation function $\rho_{2(d-k)-1}$ tends to $\frac{2k+l-d+1}{d}$ in (\ref{teh53}).  Therefore $\frac{J_2}{a^{2(d-k)-1}}$ tends to  
\begin{align}
\nonumber I^{\prime}_{k} \frac{(d-1)(d-k)!^{2}}{d^{2(d-k)-1}} \sum_{l = 0}^{d-k-1} \binom{d-2}{l,d-k-l-1,d-k-l-1,2k+l-d} = \\
\nonumber I^{\prime}_{k} \frac{(d-1)(d-k)!^{2}}{d^{2(d-k)-1}} \binom{d-2}{d-k-1}^{2},
\end{align}
where the last equation holds because both sides represent number of possibilities how to choose two sets of distinct $d-k-1$ elements from total $d-2$ elements.
This corresponds to $\lim_{a\rightarrow\infty}\frac{1}{a^{2(d-k)-1}} \var G_{d-k}(\mu_{a}^{(d)})$ in (\ref{teh6}).
\hfill $\Box $
\begin{remark} The standardization of the variance in (\ref{teh6}) corresponds to the Poisson case in Subsection 5.1. It remains to investigate whether the central limit theorem also holds in our case. \end{remark}

\noindent {\bf Acknowledgement}

\noindent This work was supported by the Czech Science Foundation, grant P201-10-0472, by grant SVV 260225 of Charles University in Prague and
 by the project 7AMB14DE006, funded by the Czech Ministery of Education and the German Academic Exchange Service.

\medskip

\begin{thebibliography}{99}
\footnotesize
\bibitem{RefA}
{\sc Baddeley A.} (2007). Spatial point processes and their applications, Stochastic geometry, Lecture Notes in Math., vol. 1892, Springer, Berlin, 1-­75. 
\bibitem{RefB}
{\sc Bene\v{s}, V. and Zikmundov\'{a}, M.} (2014). Functionals of spatial point processes having a density with respect to the Poisson process, {\em Kybernetika} {\bf 50,} 896--913.
\bibitem{RefD}
{\sc Decreusefond, L. and Flint, I.} (2014). Moment formulae for general point processes, arXiv:1211.4811v1 [math.PR]; {\em C. R. Acad. Sci. Paris, Ser. I} {\bf 352,} 357--361. 
\bibitem{RefG}
{\sc Georgii, H.O. and Yoo, H.J.} (2005). Conditional intensity and Gibbsianness of determinantal point processes, {\em J. Statist. Phys.} {\bf 118,} 55--83.
\bibitem{RefGM}
{\sc Geyer, C. J. and M\o ller J.} (1994). Simulation and likelihood inference for spatial point processes, {\em Scandinavian Journal of Statistics}, {\bf 21,} 359--373.
\bibitem{RefK}
{\sc Kendall, W. S., van Lieshout, M. and Baddeley, A.} (1999). Quermass-interaction processes: conditions for stability. {\em Adv. Appl. Prob.} {\bf 31,} 315--342.
\bibitem{RefJ}
{\sc Last, G. and Penrose, M.D.} (2011). Poisson process Fock space representation, chaos expansion and covariance inequalities, {\em Probab. Th. Relat. Fields} {\bf 150,} 663--690. 
\bibitem{RefL}
{\sc Last, G., Penrose, M.D., Schulte, M. and Th\"ale, C.} (2014). Moments and central limit theorems for some multivariate Poisson functionals, {\em Adv. Appl. Probab.} {\bf 46,} 348--364. 
\bibitem{RefM}
{\sc M\o ller, J. and Helisov\'{a}, K.} (2008). Power diagrams and interaction processes for unions of discs, {\em Adv. Appl. Probab.} {\bf 40,} 321--347.
\bibitem{RefP}
{\sc Pawlas, Z.} (2014). Self-crossing points of a line segment process, {\em Methodol. Comp. Appl. Probab.} {\bf 16,} 295--309. 
\bibitem{RefT}
{\sc Peccati, G. and Taqqu, M.S.} (2011). Wiener chaos: Moments, Cumulants and Diagrams, Bocconi Univ. Press, Springer, Milan. 
\bibitem{RefR}
{\sc Reitzner, M. and Schulte, M.} (2013). Central limit theorems for $U$-statistics of Poisson point processes, {\em Annals of Probab.} {\bf 41,} 3879--3909.
\bibitem{RefS}
{\sc Schneider, R. and Weil, W.} (2008). Stochastic and Integral Geometry, Springer, Berlin.
\bibitem{RefX}
{\sc Schreiber, T. and Yukich, J.} (2013). Limit theorems for geometric functionals of Gibbs point processes, {\em Ann. de l'Inst. Henri Poincar\'{e} - Probab. et Statist.}
{\bf 49,} 1158--1182.
\end{thebibliography}


\end{document}